\documentclass[cup7b]{cupbook}
\usepackage{datpatch}
\usepackage{amsmath}

\usepackage{graphics,lscape}
\DeclareFontEncoding{OT2}{}{} 
  \newcommand{\textcyr}[1]{%
    {\fontencoding{OT2}\fontfamily{wncyr}\fontseries{m}\fontshape{n}%
     \selectfont #1}}
\newcommand{\Sha}{{\mbox{\textcyr{Sh}}}}

\begin{document}
\def\Q{{\bf Q}}

\author[J. Brian Conrey {\it et al.}]
{J. Brian Conrey \\ American Institute of Mathematics
and University of Bristol
\and
Michael O. Rubinstein \\ University of Waterloo
\and Nina C. Snaith \and Mark Watkins \\ University of Bristol}

\chapter{Discretisation for odd quadratic twists}

\vspace{24pt}
\begin{abstract}
The discretisation problem for even quadratic twists is almost understood,
with the main question now being how the arithmetic Delaunay heuristic
interacts with the analytic random matrix theory prediction. The
situation for odd quadratic twists is much more mysterious, as the height
of a point enters the picture, which does not necessarily take
integral values (as does the order of the Shafarevich-Tate group).
We discuss a couple of models and present data on this question.
\end{abstract}

\section{Introduction}
Let $E: y^2=x^3+Ax+B$ be a fixed rational elliptic curve, and consider
the sets $S^{+}(X)$ and $S^{-}(X)$ of quadratic twists of $E$ that
contain respectively the even\footnote[1]{
 A twist is even if the order of vanishing of its $L$-function at $s=1$
 (that is, its analytic rank) is even, which is the same as saying that
 the sign of its functional equation is~$+1$; similarly for odd twists.}
and odd twists $E_d : dy^2=x^3+Ax+B$ with $|d|<X$ a fundamental discriminant.
For even twists, the Birch--Swinnerton-Dyer conjecture \cite{BSD}
states that $$L(E_d,1)={\Omega_d} {g_d\cdot \#\Sha_d\over T_d^2}$$
where $\Omega_d$ is the real period, $g_d$ is the global Tamagawa number,
$\Sha_d$ is the Shafarevich-Tate group,\footnote[2]{
 We allow the order to be zero,
 in which case we suspect a curve of higher rank.}
and $T_d$ is the order of the torsion subgroup,
all of these quantities being with respect to the quadratic twist~$E_d$.
Random matrix theory applied with orthogonal symmetry \cite{ckrs} predicts that
\begin{equation}
{\rm Prob}\bigl[L(E_d,1)\le x\bigr] \approx x^{1/2}(\log x)^{3/8}
\quad\text{as}\>\> x\rightarrow 0,
\label{eqn:prob0}
\end{equation}
where we use the $\approx$ notation to indicate
that the quotient of the two sides tends
to an unspecified constant that depends on~$E$.
Since $\#\Sha_d$ is a square while $g_d$ and $T_d$
are well-understood integers,
we get a discretisation from \eqref{eqn:prob0} --- we expect
that $L(E_d,1)=0$ if, say, we have that~$L(E_d,1)\le g_d\Omega_d/T_d^2$.
Because $\Omega_d$ essentially acts like $\approx 1/\sqrt{|d|}$,
this gives a rough prediction that
$${\rm Prob}\bigl[L(E_d,1)=0\bigr] \approx (\log |d|)^{C}/|d|^{1/4}$$
as~$|d|\rightarrow\infty$, where the constant $C$ is well-understood,
largely dependent on the rational $2$-torsion structure of~$E$.
Finally, these heuristics lead to a conjecture
about the number of positive rank twists in $S^{+}(X)$, namely that
there should be about $\approx X^{3/4}(\log X)^{C}$
of them as~$X\rightarrow\infty$.

The situation is somewhat different for odd twists; here we have
that $L(E_d,1)=0$ from the functional equation, and now
the BSD conjecture takes into account the regulator~$R_d$:
$$L'(E_d,1)={\Omega_d} {g_d\cdot R_d\#\Sha_d\over T_d^2}.$$
This regulator is rather mysterious, and, as in the case of
regulators and class numbers for real quadratic fields,
does not seem totally disjoint from the Shafarevich-Tate group.
The heuristic of Delaunay \cite{delaunay} gives some idea of how
we might expect $\#\Sha$ to be distributed, but for the regulator
we have only the lower bound of size $c\log |d|$
of Silverman \cite{silverman} and a conjectured
upper bound\footnote[3]{
 Assuming BSD and GRH we essentially get Lang's conjecture;
 in place of GRH, by bounding $L'(E_d,1)$ via convexity,
 we get a crude upper bound of~$|d|^{1+\epsilon}$.}
of $|d|^{1/2+\epsilon}$ of Lang~\cite{lang}.

Also, the analogue of \eqref{eqn:prob0} has a different exponent;
we have\footnote[4]{
 The exponent on the logarithm is $-r^2/2+r/2+3/8$, where $r$ is the
 order of the zero enforced at~$s=1$; see \cite{snaith} for the
 general case, and \cite{snaith1} for the case~$r=1$.}
\begin{equation}
{\rm Prob}\bigl[L'(E_d,1)\le x\bigr] \approx x^{3/2}(\log x)^{3/8}
\quad\text{as}\>\> x\rightarrow 0.
\label{eqn:prob}
\end{equation}
In analogy with the class number problem\footnote[5]{
 Note, however, that our $L$-values are at the center of the critical
 strip, while those for the class number problem are at the edge.}
we might be so bold as to guess that $R_d\#\Sha_d$ is always large
if nonzero, say as big as $|d|^{1/2-\epsilon}$.
Since $\Omega_d$ acts like~$\approx 1/\sqrt{|d|}$, this then implies
that~$L'(E_d,1)\gg 1/|d|^{\epsilon}$. More generally, we might guess that
\begin{equation}
\text{>}\quad
L'(E_d,1)\gg 1/|d|^\theta\qquad\text{for curves of analytic rank 1}\quad ?
\label{eqn:SNC}
\end{equation}
at least in a statistical sense, that is, there are disproportionately
few twists with nonzero $L'$-values smaller than this.

From this, in analogy with the above argument
(and ignoring logarithmic factors)
we obtain that as $|d|\rightarrow\infty$ we have
$${\rm Prob}\bigl[L'(E_d,1)=0\bigr] \approx 1/|d|^{3\theta/2},$$
so that the number of twists of rank greater than 1 should be
about $X^{1-3\theta/2}$ as $X\rightarrow\infty$.
We now proceed to give models and data which suggest
various values for~$\theta$. Note that the only provable
(assuming~BSD) bound is that~$L'(E_d,1)\gg 1/\sqrt{|d|}$,
which would lead to a prediction
of only $X^{1/4}$ odd twists of rank greater than 1.
However, for an infinite family of curves $E$ and under
the assumption of the Parity Conjecture,
Rubin and Silverberg \cite[8.2]{rs} can prove that there
are $\gg X^{1/3}$ twists of rank at least~3.

The above conjecture~\eqref{eqn:SNC} implies that $R_d$ and $\Sha_d$ are
linked in a mysterious way; if we have a generator of small height
(so that~$R_d$ is small), then this tends to make $\Sha_d$ be
larger than general.
The constructions of Rubin and Silverberg by their very nature yield
points that are of height that is polynomial in $\log |d|$ --- indeed,
almost any parametrised family will have this feature, as writing down
points of larger height is not feasible.
These facts together suggest that by taking families
with small generators we can generate large values of~$\Sha$.
However, this does not work quite so simply
in practise --- we do get large values of~$\Sha$,
but not always (as we will see in Section~\ref{sec:data}).
This is one of the reasons why we might suggest a statistical version 
of~\eqref{eqn:SNC} rather than a universal lower bound.

\section{A model from Heegner points (largely due to Birch)}
Suppose that $E$ has rank zero and $d<0$ is a fundamental
discriminant that is a square modulo~$4N$, where $N$ is the
conductor of $E$, and also assume for simplicity that $\gcd(d,6N)=1$.
By work of Gross and Zagier~\cite{gz}, we have a construction
for a point $P_d$ on $E_d$ that gives a torsion point precisely
when the rank of $E_d$ is greater than 1; indeed, the height~$\lambda$
of the constructed point is proportional to $L'(E_d,1)$:
$$\lambda(P_d)={\sqrt{|d|}\over 4\Omega_{\rm vol}}L(E,1)L'(E_d,1),$$
where here $\Omega_{\rm vol}$ is the area of the fundamental
parallelogram associated to a minimal model for~$E$.
When the rank of $E_d$ is 1, the point $P_d$ has infinite order
but is not in general a generator of the free part of the group
of rational points; the index of $P_d$ depends on $\#\Sha_d$,
but cancels out in the end.

The construction of the point $P_d$ goes via class-field theory;
we get a point $U_d$ over the Hilbert class field via a complex
multiplication result largely due to Shimura, and then sum the
conjugates to get a point first in the imaginary quadratic field
$\Q(\sqrt{d})$ and then in $\Q$ itself.
The number of conjugates of $U_d$ in the Hilbert class field is essentially
the class number $h$ of $\Q(\sqrt{d})$.
These points, all being conjugate, have the same height.
To get the height of the resulting point in $\Q$,
we model the situation by assuming that we are summing $h$ unit
vectors in \hbox{$h$-dimensional} space; this leads to the prediction
that the height is almost surely close to $h$ which is of size $\sqrt{|d|}$.
If we assume that the height of $U_d$ is not too small
we then get a prediction that $L'(E_d,1)\gg 1/|d|^\epsilon$,
leading to about $X^{1-\epsilon}$ twists in $S^{-}(X)$
which have rank~3 or greater. However, it is not clear why the height
of~$U_d$ might not be of size $1/|d|^C$~itself, as its coordinates
are in a field whose degree is of size~$\sqrt{|d|}$.

We can try to test the validity of this model by taking $d$
with $L'(E_d,1)$ small and then computing the height of the point~$U_d$
in the Hilbert class field. However, when the class field has large degree
(that is, when the class number is large), it will be difficult to
recognise the coordinates of~$U_d$, so we cannot take $|d|$ too large here.
We were thus unable to generate enough examples
to perform any real test of the model.

\section{Alternative ideas}

A less profound idea is to assert that the connection between
rank 1 and rank 3 twists should be the same as the connection between
rank~0 and rank 2 twists, at least to first approximation.  Heuristics and
random matrix theory \cite{ckrs} give $X^{3/4+\epsilon}$ rank 2 curves
amongst even quadratic twists up to $X$.
If we thus guess that there about $X^{3/4}$ twists of rank 3 up to~$X$,
via reverse-engineering the argument of two sections previous,
this could then be used to determine a value of $\theta=1/6$.

We note that there are two random matrix models that have been proposed for
modeling the zeros of $L$-functions associated with elliptic curves.
The prediction \eqref{eqn:prob} of Snaith \cite{snaith}
is extended to higher ranks by looking
at a zero-dimensional subset of $SO(\text{even})$ (for even twists)
or $SO(\text{odd})$ (for odd twists) with $r$ eigenvalues
conditioned to lie at~1.
This model predicts $\text{Prob} [L^{(r)}(E_d,1)\leq x] \approx
x^{r+1/2}(\log x)^{-r^2/2+r/2+3/8}$.  In contrast, Miller \cite{miller2}
has proposed his Independent Model, with eigenvalue distribution
decomposing as a sum of $(2 \lfloor r/2 \rfloor+1)$ point-masses and
the eigenvalue distribution of the symmetry group $SO(\text{even})$
or~$SO(\text{odd})$.
In this case the $r$th derivative analogue of \eqref{eqn:prob0}
and \eqref{eqn:prob} is given by \eqref{eqn:prob0} for $SO(\text{even})$
symmetry and \eqref{eqn:prob} for $SO(\text{odd})$ symmetry.
There is both theoretical evidence \cite{miller,young} and numerical
data~\cite{miller2} that the 1- and 2-level densities of zeros follow
Miller's Independent Model for $L$-functions associated with parameterised
families of elliptic curves with $r$ constructed points that generate
the infinite part of the Mordell-Weil group.
But there is no evidence to suggest that the Miller model
should hold in the case of quadratic twists, and in fact the
exponent~3/2 in~\eqref{eqn:prob} is supported by the shape
of the value distribution of~$L'(E_d,1)$ (see Figure~\ref{fig:graph})
as well as by the results in Section~\ref{sec:qtwist}.
This illustrates that for odd twists the zero of $L(E_d,s)$ at $s=1$ is
apparently not independent --- in contrast to a case of Young's \cite{young}
where the zero was the result of a constructed rational point on the
elliptic curve.

Finally there is a model due to A.~Granville.
Let $E$ be a fixed elliptic curve given by the model~$y^2=x^3+Ax+B$.
Here we make a heuristic for the number of integral points
$(d,u,v,w)$ with $dw^2=v(u^3+Auv^2+Bv^3)$ and~$D<|d|<2D$ and $X<|u|,|v|<2X$.
There are about $\approx X^2$ such $(u,v)$-pairs,
and each leads to a right-hand side which is of size~$X^4$.
The number of integers that are of size $X^4$
and are $d$ times a square with $D<|d|<2D$ is~$\approx D\sqrt{X^4/D}$,
and thus the probability that an integer of size~$X^4$
is of this form is $\approx \sqrt{DX^4}/X^4$.
Multiplying this by our~\hbox{$\approx X^2$}
possibilities for~$(u,v)$, we get a total of~\hbox{$\approx \sqrt D$}
integral solutions, independent of~$X$.
Summing this dyadically over~$X$, we get~\hbox{$\approx\sqrt D\log Y$}
total solutions up to~$Y$, and switching to logarithmic heights,
we get that the
number of points of height less than~$H$ on the $D$ twists of~$E$
is~$\approx H\sqrt D$. We then note (under GRH) that $E_d$ has
regulator at most size~$|d|^{1/2+\epsilon}$;
if $E_d$ is of rank~3, since a random $3$-dimensional lattice
of this covolume should have a vector whose length is of
size~$(|d|^{1/2+\epsilon})^{1/3}$, we then expect a point of
height less than $|d|^{1/6+\epsilon}$ on~$E_d$.
From the above, we expect no
more than about $X^{1/2+1/6+\epsilon}$ such twists up to~$X$.
The prediction of~$\approx H\sqrt D$ such $(d,u,v,w)$-tuples can be proved
via a sieve argument for small~$H$, but is more dubious for large~$H$.
Indeed, with just one twist of rank~$r$ with generator of
maximal height~$h$, we get $(H/h)^{r/2}$ points of height less than~$H$;
with $r=3$ and $H\rightarrow\infty$ we outdo the linear growth
predicted by the model. However, we only need $H$ to be a small
power of~$D$, and it is unclear how far the heuristic can be pushed.
Note that the obvious generalisation of this heuristic
predicts an upper bound of $X^{1/2+1/2r+\epsilon}$
for the number of rank $r$~twists.

\section{Data\label{sec:data}}

\begin{figure}[h]
\begin{center}
\hspace{-20pt}
\scalebox{0.65}{\includegraphics{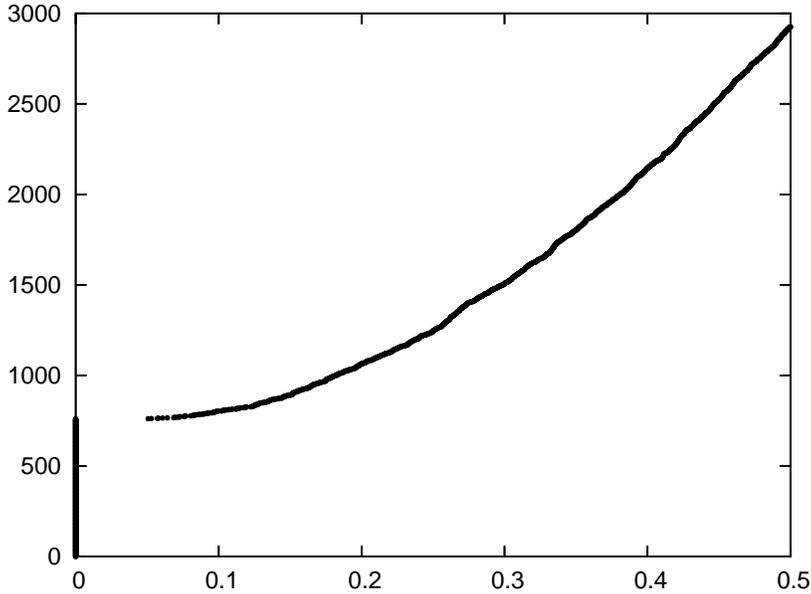}}
\vspace{-24pt}
\caption
{Cumulative $L'$-distribution for odd twists of $X_0(11)$ for $|d|<10^6$.
\label{fig:graph}}
\end{center}
\end{figure}

We now give tables and graphs that concern the above heuristics
and conjectures. In our first graph (Figure~\ref{fig:graph}),
we plot the $L'$ values for odd twists of $X_0(11)$ with~$|d|<10^6$.
We are most concerned with the behaviour as $L'\rightarrow 0$, so we
zoom in on this point; there are about 300000 total curves, of which
760 have~$L'=0$.

\begin{figure}
\begin{center}
\hspace{-20pt}
\scalebox{0.65}{\includegraphics{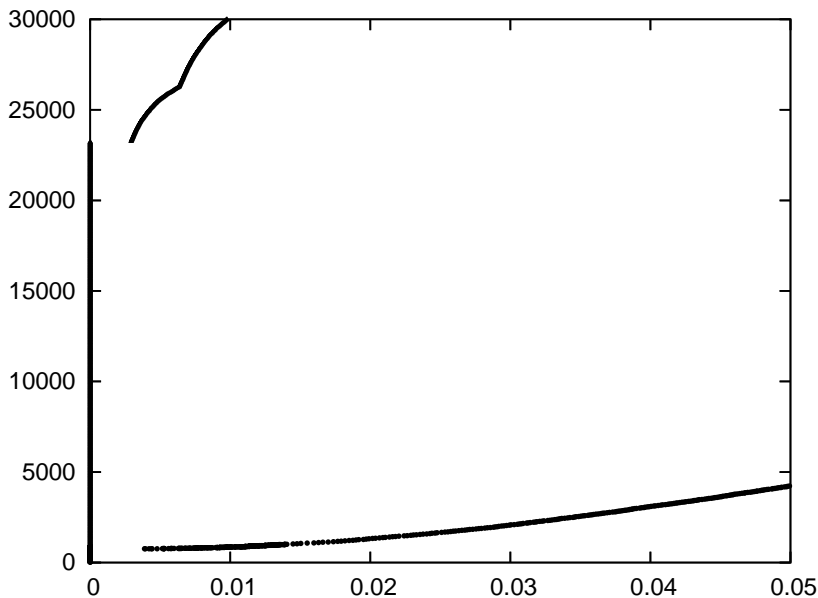}}
\vspace{-24pt}
\caption
{Cumulative distributions for $L$ (top) and normalised~$L'$ for $|d|<10^6$.
\label{fig:graph2}}
\vspace{-12pt}
\end{center}
\end{figure}

Looking at this graph, it looks as though there is an abrupt cutoff.
We find that the smallest nonzero value of $L'(E_d,1)$ is about
$0.051$ for~$d=477121$. However, it should be noted that it might be
superior to look at the distribution of~$L'(E_d,1)/\bigl(\log |d|\bigr)$,
due to the fact that the average value of $L'(E_d,1)$ is proportional
to $\log |d|$ (see~\cite{BFH,iwaniec,murty}). This changes the picture
quantitatively (see Figure~\ref{fig:graph2}),
as the gap size becomes comparable to that of the
$L$-distribution at the top of the graph.\footnote[6]{
 It can be noted that $\log |d|$ is about size $|d|^{1/6}$
 for our~$d$, and thus it becomes difficult to distinguish
 in our data between a logarithm and a power of~$d$.}

We compare the situation between even and odd twists.
For $|d|<10^6$ there are about 30 times more even twists with $L(E_d,1)=0$
than odd twists with~$L'(E_d,1)=0$; however this factor of~30
is dependent on our cutoff of~$10^6$, and as we note below,
it is not clear what happens asymptotically.
If we restricted our range of~$d$ to a shorter interval,
say $9\cdot 10^5<|d|<10^6$, then the upper graph of $L$-values would
be close to steplike, since the size of~$d$ is the only continuous
variable in the BSD formula. However, the lower graph would still be
rather smooth, since in the rank 1 case the regulator cannot be modelled
as a discrete variable.

Letting $S_0^{-}(X)$ be the subset of $S^{-}(X)$ with $L'(E_d,1)=0$,
if we believe that $\#S_0^{-}(X)\sim cX^A(\log X)^B$ we can try to
fit the data to get the exponent~$A$.
For $X_0(11)$ there are 760 odd twists with $L'=0$ with $|d|<10^6$.
The best-fit exponent for the data is~$A=0.86$, though
if we just look at the last 380 curves, we get $A=0.82$. The
computations of Elkies\footnote[7]{
 He divides even fundamental discriminants by~4,
 and so has different curve counts.} \cite{elkies}
for $X_0(32)$ go up to $10^7$, and give $A=0.84$ overall and $A=0.80$
for the last half of the data;
of course, we are ignoring log-factors, so $A=0.75$ is quite reasonable.
For $X_0(14)$ we get $A=0.94$ and for $X_0(15)$ we get $A=0.95$.
These might seem large, but Elkies has $A=0.93$ at $10^6$
before it drops significantly as indicated above.
Also, since $X_0(14)$, $X_0(15)$, and $X_0(32)$ all have
nontrivial $2$-torsion while $X_0(11)$ does not, we might expect
the exponent of the logarithm to be larger for them, which could lead
to a larger observed value of~$A$ across the range of our dataset.
For comparison with the even twist case,
the dataset of Rubinstein \cite{rubinstein-data} for the number
of rank 2 imaginary quadratic twists of $X_0(11)$ has best-fit exponents
of about $0.89,0.86,0.84$ up to $10^6,10^7,10^8$, while we expect the
exponent to be~$0.75$.

To get a dataset of twists with points of small height,
we looked at the $d$th twist of $y^2=x^3-1$ for~$d=t^3-1$;
the curve $dy^2=x^3-1$ will have the point~$(t,1)$ whose height
is of size of~$\log d$.
As mentioned above, if \eqref{eqn:SNC} holds, we would expect
such curves to have large values of~$\#\Sha_d$. Though we get some large
examples like $t=624$ and $d=242970623$ for which~$\#\Sha_d=47^2$,
this idea does not always work so well.
For instance, with $t=810$ and $d=531440999$ we have~$\#\Sha_d=1$,
where here we have~$L'(E_d,1)\approx 0.0315$;
similarly $t=902$ and $d=733870807$ has~$\#\Sha_d=1$,
though in this case $L'(E_d,1)\approx 0.0546$ is not quite so small.
Note also that the results of Delaunay and Duquesne \cite{dd}
for curves connected to the simplest cubic fields show $\#\Sha=1$
to occur quite often.

More extensive experiments using techniques similar to those of Elkies
are planned --- indeed, it would be nice to have data for the odd twists
comparable to that which \cite{ckrs2} has for even twists. Up to this point,
our experiments for odd twists have simply computed the value of $L'$ for
every twist up to $X$ and so takes $X^2$ total time, while the method
of Elkies takes $X^{3/2}$ time, as does\footnote[8]{
 With convolution techniques this can be reduced to essentially linear time,
 which is one reason why we seek to improve on \cite{elkies} via $p$-adic
 computations and $\Theta$-series.}
the computation of~\cite{ckrs2}.

\subsection{Quadratic twists in arithmetic progressions\label{sec:qtwist}}
We can note that the computations of Elkies \cite{elkies} already give
indirect evidence that \eqref{eqn:prob} is probably correct.
While Elkies notes a strange discrepancy in the counts $E_d$ with rank~3
for $d$ modulo~16, in fact, as explained in the last section of~\cite{ckrs},
we expect such discrepancies for all (prime) moduli~$p$ whose
Frobenius trace $a_p$ is nonzero.
In particular, of the $d$ with $E_d\in S_0^-(X)$
we expect that the number of nonzero quadratic residues mod~$p$
is not the same as the number of quadratic nonresidues.
The derivation in \cite{ckrs} gives a ratio of
$\bigl({p+1+a_p\over p+1-a_p})^k$ where the exponent $k=-1/2$ is taken
to be the rightmost pole of the distribution function; in the rank 1 case,
the corresponding calculation of \cite{snaith} implies that we should
take~$k=-3/2$. This is a reasonably testable prediction, given that the
dataset of Elkies has 8740 curves. In Table~\ref{table:ratio} we give
the results for some primes that are 1~mod~4; since $a_p=0$ for other
odd primes the ratio should be~1, and indeed it is always quite close.
Here the $R$ and $N$ columns count the $d$ for which $E_d$
has rank 3 and $d$ is respectively a nonzero quadratic residue and
a quadratic nonresidue mod~$p$, while the $E$ column calculates their
experimentally-determined ratio, and $C$ is the conjectured ratio from the
above with~$k=-3/2$.

\begin{table}[h]
\caption
{Effects of residuosity in arithmetic progressions
 for rank 3 quadratic twists for the congruent number curve (data from Elkies)
\label{table:ratio}}
\begin{center}
\begin{tabular}{|r|c|c|c|c|}\hline
$p$&$R$&$N$&$E$&$C$\\\hline
  5& 4240& 1951&  2.17& 2.83\\
 13& 1827& 5580&  0.33& 0.25\\
 17& 3186& 4197&  0.76& 0.72\\
 29& 5873& 2249&  2.61& 2.83\\
 37& 4451& 3820&  1.17& 1.17\\
 41& 2711& 5411&  0.50& 0.48\\
 53& 2672& 5723&  0.47& 0.45\\
 61& 5239& 3245&  1.61& 1.63\\
 73& 4696& 3688&  1.27& 1.28\\
 89& 3648& 4828&  0.76& 0.72\\
 97& 2958& 5526&  0.54& 0.57\\
929& 4836& 3876&  1.25& 1.16\\
937& 4679& 4035&  1.16& 1.13\\
941& 4807& 3922&  1.23& 1.20\\
953& 4196& 4524&  0.93& 0.92\\
977& 4791& 3929&  1.22& 1.21\\
997& 4019& 4712&  0.85& 0.83\\
\end{tabular}
\end{center}
\end{table}

Note that the fit is not as tight for small primes; indeed
this also shows up in the even rank case, even when accounting for
a secondary term as in~\cite{cprw}. Given our dataset size,
the confidence interval width for the experimental value
is about $0.1$ across most of our data range.
If we take all the primes up to 1000 and do a fit for the best~$k$, we get a
result of $-1.41$, which is reasonably close to our expected value of~$-3/2$.
This gives us a modicum of confidence that \eqref{eqn:prob} is correct;
we hope a consideration of the secondary term will give an even better fit.

\subsection{Beyond twists}
\begin{figure}[h]
\begin{center}
\hspace{-20pt}
\scalebox{0.65}{\includegraphics{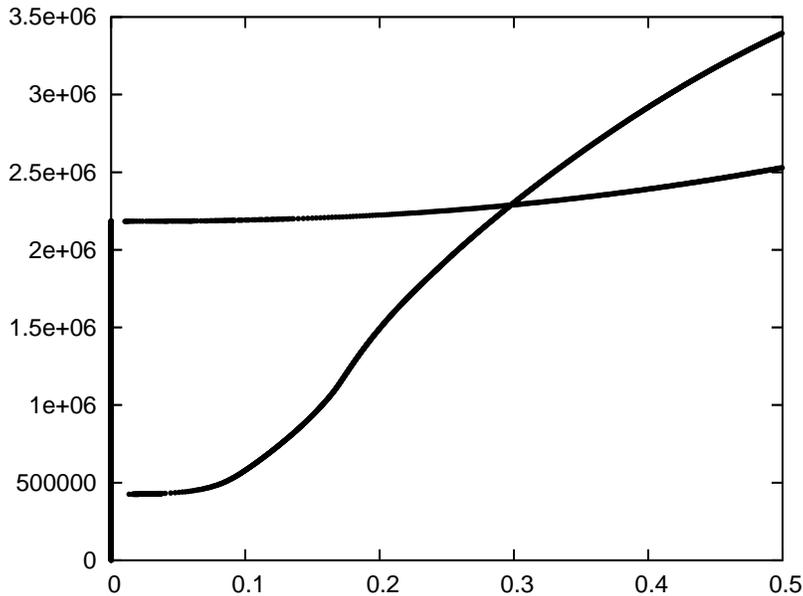}}
\vspace{-24pt}
\caption
{Cumulative distributions for $L$ and normalised-$L'$ for ECDB curves.
The plot going from the lower-left to the upper-right is
that for~$L'$.\label{fig:graph3}}
\end{center}
\end{figure}

To go further, we can look at generic elliptic curves (rather than just
twists); for this the database of Stein and Watkins \cite{sw} is useful.
Here we might guess some bound like $L'(E,1)\gg 1/|\Delta|^{\theta/6}$
in analogy with the prediction \eqref{eqn:SNC} of
$L'(E_d,1)\gg 1/|d|^{\theta}$ for quadratic twists.\footnote[9]{
 This analogy comes from the fact that the discriminant grows
 like $d^6$ in quadratic families, and our impression is that the
 discriminant is better than the conductor as a measure of the
 likelihood that the $L$-derivative vanishes. Actually we might suspect the
 real period to be the most significant datum in general,
 but it should be approximately~$|\Delta|^{1/12}$ up to log-factors.
 In any case, considering the conductor is more difficult, even with
 the ABC conjecture.}
However, as above, we really have no idea how to
generate a good value of~$\theta$.
The Stein-Watkins database (ECDB) has 11372286 curves of
prime conductor less than $10^{10}$ (we make the choice of prime
conductor so as to exclude twists from our data; looking at other
curves does not change the result too much),
of which 5253162 have analytic rank~1. The minimal $L'$-value
for these curves is about $0.193$ for the curve\footnote[10]{
 Here and below a curve $y^2+a_1xy+a_3y=x^3+a_2x^2+a_4x+a_6$
 is denoted by $[a_1,a_2,a_3,a_4,a_6]$.}
$[0,0,1,-76931443,-259719125220]$ of conductor 8519438341.
We get\footnote[11]{
 The usual caveats about not being able to prove that a curve
 actually has analytic rank~$r$ when $r\ge 4$ apply here.}
423944 curves of analytic rank~3, and 1296 of analytic rank~5.
In Figure~\ref{fig:graph3} we again see fewer curves with small
normalised $L'$-value with the normalised
gap for $L'$ about as big as that for~$L$.

It was noted to us by N.~D.~Elkies that the small values of~$L'$
correspond to curves with large cancellation between $c_4^3$ and $c_6^2$.
See Table~\ref{tbl:L1small} for the smallest values of $L'$ in the database.
For the even rank case, the smallest 85 $L$-values all come from Neumann-Setzer
\cite{neumann,setzer} curves
(with conductor of the form $u^2+64$), with the next smallest coming
from $[1,1,1,-2413424773,-45636080008772]$ of conductor 6375846313;
these thus similarly exhibit large cancellation between $c_4^3$ and $c_6^2$.
Indeed, many of the curves come from families similar to those investigated
by Delaunay and Duquesne~\cite{dd}.

\begin{table}[h]
\caption
{Small $L'$-values for prime conductor curves in the ECDB\label{tbl:L1small}}
\begin{center}
\begin{tabular}{|c|r|l|}\hline
$\>L'\>$&conductor&equation\\\hline
$0.193$&$8519438341$&$[0,0,1,-76931443,-259719125220]$\\
$0.217$&$8072290789$&$[0,-1,1,-168735150,843694875000]$\\
$0.218$&$7807742161$&$[1,0,0,-162115427,794469530026]$\\
$0.219$&$7598316169$&$[1,-1,1,-157763487,762746660718]$\\
$0.219$&$972431659$&$[1,-1,0,-42359524,-106103907983]$\\
$0.220$&$7344220789$&$[1,-1,1,-153528564,732242039802]$\\
$0.225$&$6436262197$&$[1,-1,1,-133616676,594515948970]$\\
$0.226$&$6347138731$&$[0,1,1,-131764782,582122479302]$\\
$0.226$&$2829273949$&$[1,-1,1,-119862711,-505066414494]$\\
$0.229$&$5907969559$&$[1,-1,1,-122639979,522783273972]$\\
\end{tabular}
\end{center}
\end{table}

Following a suggestion of A.~Venkatesh, we can consider whether
all the small $L'$ values (possibly including~$L'=0$) essentially come
from a small number of parametrised families.
We can make a heuristical argument against the
analogous claim that all rank 2 curves should come from parametrised families.
A heuristic of Watkins \cite{rank2} gives that there should be at least
$X^{19/24-\epsilon}$ curves of analytic rank 2 with conductor less than~$X$,
whereas we expect\footnote[12]{
 This type of heuristic appears (though not explicitly) in the work
 of Elkies and Watkins~\cite{antsVI}. They only consider small
 generators that are integral, but by passing to rationality
 we only lose logarithmic factors.}
there only to be about $X^{2/3+\epsilon}$ curves with two small generators.

We can go to curves of larger rank and look at the distribution
of $L''(E,1)/2!$ and $L'''(E,1)/3!$ for curves
of (analytic) rank 2 and 3  in the database.
If we ignore various examples of small conductor,
the smallest value of $L''(E,1)/2!$ for a curve of larger
conductor is about $1.554$ for the curve
$[0,0,1,-2664919573,-52951013063110]$ of conductor~6264757621,
where again we see the large cancellation between $c_4^3$ and~$c_6^2$.
For rank~3 the smallest value of $L'''(E,1)/3!$ for curves of larger
conductor is about $8.089$ for the curve
$[0,0,1,-7990342,8693530176]$ whose conductor is~1531408357.
Though there is large cancellation between $c_4^3$ and $c_6^2$~here,
it is not as noticeable as in the cases above;
however, the large cancellation appears again for the next-best curve
$[0,0,1,-217363231,1233466148550]$ of conductor~6352778197
for which we have~$L'''(E,1)/3!\approx 8.24$. As noted above, it is better
to divide the $L^{(r)}$-values through by the expected average value,
which is propotional to~$(\log N)^r$, before making these comparisons;
upon doing this, the listed curves of conductor~6264757621
and conductor~6352778197 have the smallest respective values.

\section{Conclusion}
Via the use of random matrix theory, we have given a link (as in the case
of rank 2 quadratic twists) between the distribution of $L'$-values and the
number of rank 3 quadratic twists, but are unable to gain much insight
into solving the discretisation problem. Although we might expect a
smooth distribution function for~$L'(E_d,1)$ (especially as it is an
analytic and not an arithmetic object), there is some evidence of a
rather abrupt cutoff in the distribution. This has led some of the authors
of this paper to conjecture \eqref{eqn:SNC} in a universal form, while others
remain more skeptical.\footnote[13]{
 It may be noted that \eqref{eqn:SNC} has been referred
 to as the ``Saturday Night Conjecture'' due its formulation 
 on a Saturday night at the Isaac Newton Institute.}
We have also discussed various methods for modelling the number of rank~3
quadratic twists of a given elliptic curve.
However, currently we do not have enough data to feel confident
in eliminating any of the suggestions.

\hyphenation{EPSRC}
\section{Acknowledgments}
This work was largely done during the programme
``Random Matrix Theory and Number Theory''
at the Isaac Newton Institute (INI), whose hospitality we enjoyed.
The first and second authors were partially supported by
a Focused Research Grant from the~National Science Foundation~(USA),
the second author by a
Natural Sciences and Engineering Research Council grant (Canada),
the third author by
a Dorothy Hodgkin Fellowship from
the Royal Society and
Engineering and Physical Sciences Research Council (EPSRC/UK) grants
GR/T00825/01 and GR/T00832/01, and 
and the fourth author by an INI Fellowship
(EPSRC grant GR/N09176/01)
and EPSRC grant GR/T00658/01.

\end{document}